\documentclass[10pt]{amsart}
\usepackage[T1]{fontenc}
\usepackage[latin1]{inputenc}
\usepackage{amsmath}
\usepackage{amssymb}
\usepackage{amsthm}
\usepackage{dsfont}
\usepackage{color}
\usepackage{enumitem}
\theoremstyle{plain}\newtheorem{theo}{Theorem}
\theoremstyle{plain}\newtheorem{cor}{Corollary}
\theoremstyle{definition}\newtheorem{rem}{Remark}
\theoremstyle{plain}
\theoremstyle{plain}
\theoremstyle{plain}
\theoremstyle{definition}

\newcommand{\N}{{\mathbb{N}}}

\begin{document}
\title[The Gibbs sampler for the 1-D Ising model]{The convergence rate of the Gibbs sampler for \\generalized $1-$D Ising model}
\author[ Amine Helali  ]{ Amine Helali} 
\address{Laboratoire d'Algèbre G\'{e}om\'{e}trie et Th\'{e}orie Spectrale LR/11/ES-53, D\'{e}partement de  Math\'{e}matique, Facult\'e des Sciences de Sfax, Universit\'{e} de Sfax, 3000 Sfax, Tunisia}
\address{Laboratiore de Math\'{e}matiques de Bretagne Atlantique UMR 6205, UFR Sciences et Techniques, Universit\'{e} de Bretagne Occidentale, 6 Avenue Le Gorgeu, CS 93837, 29238 Brest, cedex 3, France.}
\email{ amine.helali@univ-brest.fr}
\subjclass[2010]{Primary: 60J22, Secondary: 60F99, 60J10.}
\keywords{Markov chain Monte Carlo, Rate of convergence, Gibbs sampler,  Ising model.}
\begin{abstract}
The rate of convergence of the Gibbs sampler for the generalized one-dimensional Ising model is determined by the second largest eigenvalue of its transition matrix in absolute value denoted by $\beta^*$. In this paper we generalize a bound for $\beta^*$ from Shiu and Chen $(2015)$ for the one-dimensional Ising model with two states to a multiple state situation. The method is based on Diaconis and Stroock bound for reversible Markov processes. The new bound presented in this paper improves Ingrassia's $(1994)$ result.

\end{abstract}                        
\maketitle
\section*{Introduction}
The Ising model is a crude model for ferromagnetism. It is the simplest model of statistical mechanics and it has been applied in many other fields like chemistry, molecular biology and image analysis.   
The distribution of the one-dimensional Ising model with three states is:
 $$\pi(x)=\frac{1}{Z_T}\exp\Big\{ \tfrac{1}{T}\displaystyle\sum_{k=1}^{n-1} \big{(} \mathds{1}_{\{x_k = x_{k+1}\}} - \mathds{1}_{\{x_k \neq x_{k+1}\}}   ) \Big\}  \,\,\,\,\,\,\,\,\, \forall \,x = (x_1, \cdots, x_n \big{)}  \in  \chi\,\,, $$
where $\chi=\{c^{(1)}, c^{(2)}, c^{(3)}\}^n$ is the state space, $T$ is the temperature and $$Z_T =  \displaystyle\sum_{x  \in  \chi} \exp \Big\{ \tfrac{1}{T}   \displaystyle\sum_{k=1}^{n-1} \big{(} \mathds{1}_{\{x_k = x_{k+1}\}} - \mathds{1}_{\{x_k \neq x_{k+1}\}} \big{)} \Big\}$$ is the normalizing constant. Monte Carlo Markov chain $MCMC$ method is a very useful technique to draw samples from the Ising model.
Suppose that the transition probability $P(x, y)$ for an irreducible Markov chain  has $\pi$ as its invariant measure. The pair $(P, \pi)$ is said to be reversible if it verifies the detailed balance equation: 
$$Q(x,y) = \pi(x) P(x,y) = \pi(y) P(y,x)  = Q(y,x) \,\, \,\,\,\,\,\,\,\, \forall   \,x,y\in \chi . $$ 
The Gibbs sampler introduced by Geman and Geman and the Metropolis-Hastings algorithm introduced by  Metropolis et al. and Hastings (see \cite{Geman and Geman},  \cite{Metropolis et al}  and \cite{Hastings}) are the most popular Monte Carlo Markov chain methods. The matrix $P$ satisfies detailed balance and thus is symmetric with respect to the scalar product introduced by the measure $\pi$. Therefore its eigenvalues are real and can be arranged as follows:
$$ 1= \beta_0 > \beta_1 \geq \beta_2 \geq \cdots \geq \beta_{| \chi | -1} > -1.$$ 
Let $\beta^* = \max\{\beta_1, | \beta_{| \chi | -1} | \}.$ By using the total variance distance, the second largest eigenvalue in absolute value determines the convergence rate of the Markov chain. Ingrassia  gives a lower bound for $\beta_{| \chi | -1}$ and an upper  bound for $ \beta_1$ (see \cite{Ingrassia}).  \\
 This paper deals with the Gibbs sampler for the one-dimensional Ising model with multiple states. It chooses a random coordinate which is updated according to the conditional probability given the other coordinates. The resulting Markov chain is reversible and the associated transition matrix has the form:
$$P(x, y)=  \left\{
    \begin{array}{ll}
     \frac{ 1 }{ n }  \pi(y_i \mid x) \qquad \qquad \qquad \qquad \qquad \qquad \mbox{ if}\,\,  x_j \neq y_j \,\, \mbox{for\,\,all}  \,\, j \neq i
\\      
1- \frac{ 1 }{ n }  \displaystyle\sum_{i=1}^{n} \displaystyle\sum_{y_i\in  \{ c^{(1)}, c^{(2)}, c^{(3)} \}} \pi(y_i \mid x) \,\,\,\,\,\,\,\, \,\, \,\, \mbox{ if }\,\,x=y
\\
0 \qquad \qquad \qquad \qquad \qquad \qquad \qquad \qquad \,\,\,  \mbox{else}   
    \end{array}
\right.$$
where 
$$ \pi(y_i | x) =  \frac{  \pi(x_1,  \cdots , x_{i-1}, y_i,  x_{i+1}, \cdots, x_n  )}{ \displaystyle\sum_{l=1}^{3}  \pi(x_1,  \cdots , x_{i-1}, c^{(l)},  x_{i+1}, \cdots, x_n  )} .$$
Diaconis and Stroock (see \cite{Diacons Stroock}) give  a bound for the total variation distance to equilibrium in terms of $\beta^*$. We recall this result in the following theorem:
\begin{theo} [\mbox{Diaconis and Stroock 1991}]
If $P$ is a reversible Markov chain with unique invariant measure $\pi$ and $P$ is irreducible then for all $x \in \chi$ and $k  \in \N$:
$$ 4 \mid\mid P^k(x, .)  - \pi  \mid\mid_{var}^2 = \Big{(}   \displaystyle\sum_{y \in \chi } \mid P^k(x, y) - \pi(y) \mid \Big{)}^2 \leq \frac{1-\pi(x)}{\pi(x)} (\beta^*)^{2k}. $$
\end{theo}
Moreover, Diaconis and Stroock (see \cite{Diacons Stroock}) develop a method to calculate an upper bound for the second largest eigenvalue $\beta_1$  for a reversible Markov chain using geometric quantities such as the maximum degree, diameter and covering number of the associated graph. Consider the graph $\mathcal{G}(P)=(\chi, E)$ where $\chi$ is the vertex set and $E=\{(x, y)\mid P(x, y) >0\}$ is the edge set. For each pair of distinct points $x, y \in \chi$ we choose a path $\gamma_{xy}$ from $x$ to $y$, such that each edge appears at most once in a given path. The fact that $P$ is irreducible guarantees that such paths exist. Let  $\Gamma$ be the collection of all such paths $\gamma_{xy}$ (one for each pair). The geometric bound given by Diaconis and Stroock  is, 
\begin{eqnarray}\label{Geometric bound}
\beta_1 &\leq& 1 - \frac{1}{\kappa}
\end{eqnarray}
with
\begin{align}
\kappa = \underset{e\in E}\max \,\,\, Q(e)^{-1} \sum_{\gamma_{xy} \ni e} |\gamma_{xy}|\pi(x) \pi(y) \label{kappa}
\end{align}
where  $| \gamma_{xy} | $ designates the length of the path $\gamma_{xy}$ (see \cite{Diacons Stroock}).\\
In their paper Shiu and Chen (see \cite{Shiu chen}) present a method to explicitly compute the bound of Diaconis and Stroock (see \cite{Diacons Stroock}) for the Gibbs sampler for a two state one-dimensional Ising model.\\ 
In this paper we generalize the result of Shiu and Chen  to the case of the one-dimensional Ising model with three and more states (see \cite{Shiu chen}).\\
Our method is based on the idea from  \cite{Shiu chen} which consists of defining suitable paths $\gamma_{xy}$ linking each pair $(x, y)$ from the state space $\chi$ and  then to explicitly compute $\kappa$ defined in equation (\ref{kappa}) with some suitable symmetry argument. In the discussion section of the paper we compare  our bound to results from the literature. It turns out that the result generalizes the bound given in \cite{Shiu chen} to the case of the Ising model with three states (see Theorem $2$) and also to multiple states (see Theorem $3$). It also improves the bound presented by Ingrassia in \cite{Ingrassia}.
\section{Main result}
\subsection{Selection of paths}
 To be able to use the result of Diaconis and Stroock (see \cite{Diacons Stroock}) and to calculate the bound of the second largest eigenvalue, we have to fix a collection of paths connecting any configuration $x \in \chi$ to any configuration $y \in \chi$ . To get a small upper bound for $\beta_1$, we seek a small value for $\kappa$ and we should therefore use short paths $\gamma_{xy}$ to link $x$ with $y$. Moreover, we have to keep the number of paths passing through a given edge low.\\
For a pair of distinct configurations $x, y \in \chi$  there exist some increasing sequence $d_1, \cdots, d_m$ such that $x_i\neq y_i$ for  $i \in \{ d_1, \cdots, d_m \}$ and    $x_i = y_i$ otherwise.
In the same way as  Shiu and Chen (see \cite{Shiu chen})   we define a path linking a given pair $(x, y)$ as follows:
\begin{align}
(x_1, \cdots, x_n) &= (y_1, \cdots, y_{d_1-1}, x_{d_1}, x_{d_1+1}, \cdots, x_{d_2-1}, x_{d_2}, x_{d_2+1}, \cdots, x_n) \nonumber \\
&\rightarrow  (y_1, \cdots, y_{d_1-1}, y_{d_1}, x_{d_1+1}, \cdots, x_{d_2-1}, x_{d_2}, x_{d_2+1}, \cdots, x_n) \nonumber \\
&= (y_1, \cdots, y_{d_1-1}, y_{d_1}, y_{d_1+1}, \cdots, y_{d_2-1}, x_{d_2}, x_{d_2+1}, \cdots, x_n) \nonumber \\
&\rightarrow  (y_1, \cdots, y_{d_1-1}, y_{d_1}, y_{d_1+1}, \cdots, y_{d_2-1}, y_{d_2}, x_{d_2+1}, \cdots, x_n) \nonumber \\
&\vdots  \nonumber \\
& \rightarrow (y_1, \cdots, y_n).    \nonumber
\end{align}
We turn now to give an upper bound for the value of $\kappa$ defined in equation (\ref{kappa}).
\subsection{Geometric bound of the second largest eigenvalue}     
In what follows we will essentially follow the arguments from Shiu and Chen (see \cite{Shiu chen}) to find an upper bound for $$ \kappa = \underset{e \in E}\max \,\ Q(e)^{-1} \sum_{\gamma_{xy} \ni e} | \gamma_{xy} |  \pi(x) \pi(y).$$ 
Let $e=(e^-, e^+)$ be some edge from $E$ where $e^-$ and $e^+$ are two configurations from $\chi$ which differ by only one coordinate. Without loss of generality we consider the case where the site $i$  choosen to be updated passes from the color $c^{(1)}$ to the color $c^{(2)}$. Then the configurations $e^-$ and $e^+$ must have the following form: $$e^-=(z_1, z_2, \cdots, z_{i-1}, c^{(1)}, z_{i+1}, \cdots, z_n)\,\, \,\,\, \mbox{and}\,\,\, \,\,e^+=(z_1, z_2, \cdots, z_{i-1}, c^{(2)}, z_{i+1}, \cdots, z_n).$$ The transition probability for a transition from $e^-$ to $e^+$ can then be computed. A short computation shows for $i=1$:
\begin{align}
P( e^-, e^+) &= \frac{1}{n} \frac{\pi(c^{(2)}, z_2, \cdots, z_n)}{\displaystyle \sum_{j=1}^{3} \pi(c^{(j)}, z_2, \cdots, z_n)} \nonumber \\
&= \frac{1}{n} \frac{\exp\Big\{ {\tfrac{1}{T} \big{(}\mathds{1}_{\{c^{(2)} = z_{2}\}} - \mathds{1}_{\{c^{(2)} \neq z_{2}\}}  + \displaystyle \sum_{k=2}^{n-1}   \mathds{1}_{\{z_k = z_{k+1}\}} - \mathds{1}_{\{z_k \neq z_{k+1}\}}\big{)}  }  \Big\} }{  \displaystyle \sum_{j=1}^{3}  \exp{ \Big\{ \tfrac{1}{T} \big{(}  \mathds{1}_{\{c^{(j)}= z_{2}\}} - \mathds{1}_{\{c^{(j)} \neq z_{2}\}}\big{)}  }  + \displaystyle \sum_{k=2}^{n-1}   \mathds{1}_{\{z_k = z_{k+1}\}} - \mathds{1}_{\{z_k \neq z_{k+1}\}}\big{)}    \Big\}   } \nonumber \\
&=\frac{1}{n} \frac{ \exp\Big\{ {\tfrac{1}{T} \big{(}\mathds{1}_{\{c^{(2)} = z_{2}\}} - \mathds{1}_{\{c^{(2)} \neq z_{2}\}}\big{)}  }  \Big\}     }{  \displaystyle \sum_{j=1}^{3}  \exp{ \Big\{ \tfrac{1}{T} \big{(}  \mathds{1}_{\{c^{(j)}= z_{2}\}} - \mathds{1}_{\{c^{(j)} \neq z_{2}\}}\big{)}  }  \Big\}    } , \label{transition1} 
\end{align}
similarly for $i=n$:
\begin{align}
P( e^-, e^+) &=\frac{1}{n} \frac{ \exp\Big\{ {\tfrac{1}{T} \big{(}\mathds{1}_{\{ z_{n-1} = c^{(2)}\}} - \mathds{1}_{\{ z_{n-1} \neq c^{(2)}\}}\big{)}  }  \Big\}     }{  \displaystyle \sum_{j=1}^{3}  \exp{ \Big\{ \tfrac{1}{T} \big{(}  \mathds{1}_{\{ z_{n-1}=c^{(j)}\}} - \mathds{1}_{\{ z_{n-1} \neq c^{(j)}\}}\big{)}  }  \Big\}    }  \label{transition3} 
\end{align}
and for $i\in \{2, \cdots, n-1\}$:
\begin{align}
P( e^-, e^+) &=\frac{1}{n} \frac{ \exp\Big\{ {\tfrac{1}{T} \big{(}\mathds{1}_{\{z_{i-1} = c^{(2)}\}} - \mathds{1}_{\{z_{i-1} \neq c^{(2)}\}} + \mathds{1}_{\{c^{(2)} = z_{i+1}\}} - \mathds{1}_{\{c^{(2)} \neq z_{i+1}\}}\big{)}  }  \Big\}     }{  \displaystyle \sum_{j=1}^{3}  \exp{ \Big\{ \tfrac{1}{T} \big{(}\mathds{1}_{\{z_{i-1} = c^{(j)}\}} - \mathds{1}_{\{z_{i-1} \neq c^{(j)}\}} + \mathds{1}_{\{c^{(j)}= z_{i+1}\}} - \mathds{1}_{\{c^{(j)} \neq z_{i+1}\}}\big{)}  }  \Big\}    } . \label{transition2} 
\end{align}
Then we turn to compute an upper bound of $\kappa$ defined in equation (\ref{kappa}) for each class of edges. The main conclusion of this paper is given in the following theorem:  
\begin{theo}
The second largest eigenvalue eigenvalue of the Gibbs sampler for the one-dimensional Ising model with three states satisfies :
\begin{align} 
 \beta_1 < 1- 3 \times n^{-2}\frac{e^{-\frac{4}{T} }}{ 1+2e^{-\frac{4}{T} }} \label{resultat1}
\end{align}
\end{theo}
The proof of this theorem is given in section $3$. The above theorem can be generalized to the case of multiple colors  where the state space is $\chi= \{c^{(1)}, \cdots, c^{(N)} \}^n$ as follows:
\begin{theo}
The second largest eigenvalue of the Gibbs sampler for the one-dimensional Ising model with $N$ states satisfies:
\begin{align} 
 \beta_1 &< 1- N \times n^{-2}\frac{e^{-\frac{4}{T} }}{ 1+(N-1)e^{-\frac{4}{T} }} . \label{resultat2}
\end{align} 
\end{theo}
We give a sketch of the proof of this theorem in section $3$.\\
To be able to quantify the convergence rate with Theorem $1$ given  by Diaconis and Stroock (see \cite{Diacons Stroock}) we must control the smallest eigenvalue in order to bound the second largest eigenvalue in absolute value. This question is addressed in the following subsection:
\subsection{Bound for the absolute value of the second largest eigenvalue}
A theorem proved by Ingrassia (see \cite{Ingrassia}, Theorem $5.3$) gives the following lower bound  for the smallest eigenvalue:
$$ \beta_{| \chi | - 1} \geq -1 + \frac{2}{1+(C-1)e^{\frac{\Delta}{T}}}.$$
For the one-dimensional Ising model with three states, $C=3$ and $\Delta=2$. This yields  for any natural number $n> 3 / \sqrt{2}$
\begin{align}
| \beta_{| \chi | - 1} | &\leq |  -1 + \frac{2}{1+ 2e^{\frac{2}{T}}} | = 1- \frac{2}{1+ 2e^{\frac{2}{T}}} < 1- \frac{2}{3}e^{-\frac{2}{T}}  <   1- 3 n^{-2}\frac{ e^{-\frac{4}{T} }}{ 2e^{-\frac{4}{T} }+1}.\nonumber
\end{align}
In the general case where $\chi=\{c^{(1)}, \cdots, c^{(N)}\}^n$,  the parameter $C$ is equal to $N$ and Ingrassia's bound behaves for any natural number  $n> N / \sqrt{2} $  as follows: 
\begin{align}
| \beta_{| \chi | - 1} | &\leq |  -1 + \frac{2}{1+ (N-1)e^{\frac{2}{T}}} | = 1- \frac{2}{1+ (N-1)e^{\frac{2}{T}}} < 1- \frac{2}{N}e^{-\frac{2}{T}}  <   1- N \times n^{-2}\frac{e^{-\frac{4}{T} }}{ 1+(N-1)e^{-\frac{4}{T} }}.\nonumber
\end{align}
The previous considerations prove the following corollary: 
\begin{cor}
The upper bounds for $\beta_1$ given in theorem $2$ and theorem $3$ are also upper bounds for the absolute value of all eigenvalues $\{\beta_1, \cdots, \beta_{|\chi|-1}\}$ of the Gibbs sampler for the one-dimensional Ising model with three states and more respectively. 
\end{cor}
\section{Discussion} 
Ingrassia (see \cite{Ingrassia}) gives the following upper bound for the second largest eigenvalue of the Gibbs sampler: 
$$  \beta_1 \leq 1- \frac{   Z_T   }{b_{\Gamma} \,\, \gamma_{\Gamma} \,\, C \,\,  |S|  } e^{-\frac{m}{T}}.$$
In this expression $Z_T$ is the normalizing constant, $S$ is the lattice of sites, $\Gamma$ is the collection of paths,  $\gamma_{\Gamma}$ is the maximum length of each path $\gamma_{xy} \in \Gamma$,  $b_{\Gamma}$ is the maximum number of paths containing any edge of $\Gamma$, $C$ is the number of configurations that differ by only one site and $m$ is the least total elevation gain of the Hamiltonian function in the sense as described by Holley and Stroock  (see \cite{Holley and Stroock}).\\ 
In our case, we have:   $|S| = n$, $\gamma_{\Gamma}=n $, $b_{\Gamma}=3^{n-1}$, $C=3$, $Z_T \leq 3 ( 1+ 2e^{-\frac{1}{2T}} )^{n-1}$ and  $m=2$. It gives that:
$$ \beta_1 \leq 1-  n^{-2}\left(\frac{1+2e^{\frac{-1}{2T}}  }{3 } \right)^{n-1} e^{\frac{-2}{T}}.$$
This upper bound differs from the result introduced in Theorem $2$ by the multiplicative factor 
 \begin{align} \theta &= \frac{e^{\frac{2}{T}} +2e^{-\frac{2}{T}}}{3}  \left(\frac{1+2e^{-\frac{1}{2T}}  }{3 } \right)^{n-1}\nonumber 
 \end{align}
 To get improvement we need that $\theta <1$ which means $  \frac{e^{\frac{2}{T}} +2e^{-\frac{2}{T}}}{3}  \left(\frac{1+2e^{-\frac{1}{2T}}  }{3 } \right)^{n-1}<1.$\\
An elementary computation leads to: $n>  \frac{ \log \left( \frac{\exp{\left(\frac{2}{T}\right) } +2 \exp{\left(\frac{-2}{T} \right)}}{ 3}\right)}{ \log \left(\frac{3}{1+2\exp{\left(\frac{-1}{2T}\right)}  }\right)}+1.(*)$\\
For a choice of temperature $T$ near to zero we can find an integer $n$ sufficiently large which verifies $(*)$ (it is natural in the case of the Gibbs sampler where $n\sim 10^{23})$.
\begin{rem} 
The application of Ingrassia's bound to the one-dimensional Ising model with multiple states gives: $$ \beta_1 \leq 1-  n^{-2}\left(\frac{1+(N-1)e^{\frac{-1}{2T}}  }{N } \right)^{n-1} e^{\frac{-2}{T}}$$ which differs from the result introduced in Theorem $3$ by the factor $ \tilde{\theta}$ defined as follows: $$ \tilde{\theta} =  \frac{e^{\frac{2}{T}} +(N-1)e^{-\frac{2}{T}}}{N}  \left(\frac{1+(N-1)e^{-\frac{1}{2T}}  }{N } \right)^{n-1} . $$
As previous, an elementary computation leads to: $n>  \frac{ \log \left( \frac{\exp{\left(\frac{2}{T}\right) } +(N-1) \exp{\left(\frac{-2}{T} \right)}}{ N}\right)}{ \log \left(\frac{N}{1+(N-1)\exp{\left(\frac{-1}{2T}\right)}  }\right)}+1.(**)$\\
For a choice of temperature $T$ near to zero we can find an integer $n$ sufficiently large which verifies $(**)$ (it is natural in the case of the Gibbs sampler where $n\sim 10^{23})$.
\end{rem}
\section{Proofs of the main results}
\subsection{Proof of theorem $2$} We have two principle cases:\\
a) If $i \neq \{1,n\}$: According to the selection of the paths in Section $1.1$, if a path $ \gamma_{xy}$ passing through the edge $e=(e^-, e^+)$  connects $x$ with $y$, then these extremities must have the following form:
$$ x=(x_1, x_2, \cdots, x_{i-1}, c^{(1)}, z_{i+1}, \cdots, z_n) \,\,\,\,\,\,\mbox{and}\,\,\,\,\,\, y=(z_1, z_2, \cdots, z_{i-1}, c^{(2)}, y_{i+1}, \cdots, y_n).$$ 
This yields: 
\begin{align}
\pi(x)&= \frac{1}{Z_T} \exp \bigg\{  \tfrac{1}{T}\Big{(} \displaystyle \sum_{k'=1}^{i-2} (\mathds{1}_{\{x_{k'} = x_{k'+1}\}} - \mathds{1}_{\{x_{k'} \neq x_{k'+1}\}}) +  \displaystyle \sum_{k'=i+1}^{n-1} (\mathds{1}_{\{z_{k'} = z_{k'+1}\}} - \mathds{1}_{\{z_{k'} \neq z_{k'+1}\}}) \nonumber \\
&+ (\mathds{1}_{\{x_{i-1} = c^{(1)}\}} - \mathds{1}_{\{x_{i-1} \neq c^{(1)}\}} + \mathds{1}_{\{c^{(1)} = z_{i+1}\}} - \mathds{1}_{\{c^{(1)} \neq z_{i+1}\}})  \Big{)} \bigg\}. \nonumber
\end{align}
 Moreover, the probabilities $\pi(y)$ and $\pi(e^-)$ can be expressed  similarly. It follows that:
\begin{align}
Q(e)^{-1}\pi(x) \pi(y)  &= \frac{\pi(x) \pi(y)}{\pi(e^-) P(e^-, e^+)} \nonumber \\  
&=  \frac{n}{Z_T}     \bigg\{1 + \exp{ \Big\{  \tfrac{1}{T} \big{(}-\mathds{1}_{\{z_{i-1} = c^{(2)}\}} + \mathds{1}_{\{z_{i-1} \neq c^{(2)}\}} - \mathds{1}_{\{c^{(2)} = z_{i+1}\}} + \mathds{1}_{\{c^{(2)} \neq z_{i+1}\}} }\big{)} \Big\} \nonumber \\
       &\times \displaystyle\sum_{j=1,3} \exp{ \Big\{ \tfrac{1}{T} \big{(}\mathds{1}_{\{z_{i-1} = c^{(j)}\}} - \mathds{1}_{\{z_{i-1} \neq c^{(j)}\}} + \mathds{1}_{\{c^{(j)} = z_{i+1}\}} - \mathds{1}_{\{c^{(j)} \neq z_{i+1}\}}  \big{)} \Big\} }\bigg\}   \nonumber \\
       &\times \exp \Big\{ \tfrac{1}{T} \Big{(} \displaystyle \sum_{k'=1}^{i-2} (\mathds{1}_{\{x_{k'} = x_{k'+1}\}} - \mathds{1}_{\{x_{k'} \neq x_{k'+1}\}}) +   \displaystyle \sum_{k'=i+1}^{n-1} (\mathds{1}_{\{y_{k'} = y_{k'+1}\}} - \mathds{1}_{\{y_{k'} \neq y_{k'+1}\}}) \Big{)} \Big\}     \nonumber \\
       &\times \frac{  \exp{ \Big\{ \tfrac{1}{T} \big{(}\mathds{1}_{\{x_{i-1} = c^{(1)}\}} \!-\! \mathds{1}_{\{x_{i-1} \neq c^{(1)}\}}\!+\! \mathds{1}_{\{z_{i-1} = c^{(2)}\}}\!-\!\mathds{1}_{\{z_{i-1} \neq c^{(2)}\}} \!+\! \mathds{1}_{\{c^{(2)} = y_{i+1}\}}\! -\! \mathds{1}_{\{c^{(2)} \neq y_{i+1}\}} \big{)} \Big\} } }{ \exp \Big\{ \tfrac{1}{T} \big{(} \mathds{1}_{\{z_{i-1} = c^{(1)}\}}\!-\! \mathds{1}_{\{z_{i-1} \neq c^{(1)}\}}  \big{)} \Big\} } .  \nonumber
\end{align}
We introduce the notation \begin{align} (x, c^{(l)}, y) := (x_1, x_2, \cdots, x_{i-1}, c^{(l)}, y_{i+1}, \cdots, y_{n-1}, y_n) \label{notation} \end{align}  for  $l \in \{1,2,3\}$. The previous expression becomes: 
\begin{align}
Q(e)^{-1}\pi(x) \pi(y) &=n \pi(x, c^{(1)}, y) \exp{ \Big\{  \tfrac{1}{T} \big{(}-\mathds{1}_{\{c^{(1)}=y_{i+1}\}} + \mathds{1}_{\{  c^{(1)} \neq y_{i+1}\}}  + \mathds{1}_{\{z_{i-1} = c^{(2)}\}}  - \mathds{1}_{\{z_{i-1} \neq c^{(2)}\}}  }   \nonumber  \\
&+  \mathds{1}_{\{c^{(2)} = y_{i+1}\}} - \mathds{1}_{\{c^{(2)} \neq y_{i+1} \}} - \mathds{1}_{\{z_{i-1} = c^{(1)}\}} + \mathds{1}_{\{z_{i-1} \neq c^{(1)}\}} \}  \big{)} \Big\} \nonumber  \\
&+ n \pi(x, c^{(2)}, y) \exp{ \Big\{  \tfrac{1}{T} \big{(} -\mathds{1}_{\{x_{i-1} = c^{(2)}\}} + \mathds{1}_{\{x_{i-1} \neq c^{(2)}\}}  -\mathds{1}_{\{z_{i-1} = c^{(1)}\}} + \mathds{1}_{\{z_{i-1} \neq c^{(1)}\}  }    } \nonumber \\
&+ \mathds{1}_{\{x_{i-1} =c^{(1)}\}} - \mathds{1}_{\{x_{i-1} \neq c^{(1)}\}}  - \mathds{1}_{\{c^{(2)} = z_{i+1}\}} + \mathds{1}_{\{c^{(2)} \neq z_{i+1}\}} \big{)} \Big\} \displaystyle\sum_{j=1, 3} \exp \Big\{ \tfrac{1}{T} \big{(} \mathds{1}_{\{z_{i-1} = c^{(j)}\}} \nonumber \\
&- \mathds{1}_{\{z_{i-1} \neq c^{(j)}\}} + \mathds{1}_{\{c^{(j)} = z_{i+1}\}} - \mathds{1}_{\{c^{(j)} \neq z_{i+1}\}} \big{)}\Big\}.      \nonumber \\
&= n \alpha \pi(x, c^{(1)}, y) \exp \Big\{ \tfrac{1}{T} \big{(}-\mathds{1}_{\{c^{(1)}=y_{i+1}\}} + \mathds{1}_{\{  c^{(1)} \neq y_{i+1}\}}  + \mathds{1}_{\{c^{(2)} = y_{i+1}\}} - \mathds{1}_{\{c^{(2)} \neq y_{i+1} \}}  \big{)}\Big\}   \nonumber  \\
&+ n \beta \pi(x, c^{(2)}, y) \exp \Big\{ \tfrac{1}{T}\big{(} -\mathds{1}_{\{x_{i-1} = c^{(2)}\}} + \mathds{1}_{\{x_{i-1} \neq c^{(2)}\}}  + \mathds{1}_{\{x_{i-1} = c^{(1)}\}} - \mathds{1}_{\{x_{i-1} \neq c^{(1)}\}}  \big{)}\Big\} \,\,. \label{eq11}
\end{align}  
where
\begin{enumerate}[label= \roman*)]
\item  $\alpha = \exp{ \Big\{  \tfrac{1}{T} \big{(} \mathds{1}_{\{z_{i-1} = c^{(2)}\}}  - \mathds{1}_{\{z_{i-1} \neq c^{(2)}\}}  } - \mathds{1}_{\{z_{i-1} = c^{(1)}\}} + \mathds{1}_{\{z_{i-1} \neq c^{(1)}\}} \}  \big{)} \Big\}.$ \\
\item  $\beta=\exp{ \Big\{  \tfrac{1}{T} \big{(} -\mathds{1}_{\{z_{i-1} = c^{(1)}\}} + \mathds{1}_{\{z_{i-1} \neq c^{(1)}\}  } - \mathds{1}_{\{c^{(2)} = z_{i+1}\}} + \mathds{1}_{\{c^{(2)} \neq z_{i+1}\}} \big{)} \Big\} } \displaystyle\sum_{j=1, 3} \exp \Big\{ \tfrac{1}{T} \big{(} \mathds{1}_{\{z_{i-1} = c^{(j)}\}} \\- \mathds{1}_{\{z_{i-1} \neq c^{(j)}\}}+ \mathds{1}_{\{c^{(j)} = z_{i+1}\}} - \mathds{1}_{\{c^{(j)} \neq z_{i+1}\}} \big{)}\Big\}.$
\end{enumerate}
From the notation in equation (\ref{notation}) we have  $ \displaystyle  \bigcup_{(x,y):\,\,\gamma_{xy} \ni e} \Big\{  (x,c^{(1)},y), (x,c^{(2)},y), (x,c^{(3)},y)  \Big\} := \chi$. This yields
\begin{align}
Q(e)^{-1} \sum_{(x, y):\, \gamma_{xy} \ni e} |  \gamma_{xy} | \pi(x) \pi(y) &\leq \alpha n^2  \sum_{(x, y):\,\gamma_{xy} \ni e}   \pi(x, c^{(1)}, y) \exp \Big\{ \tfrac{1}{T}\big{(} -\mathds{1}_{\{c^{(1)}=y_{i+1}\}} + \mathds{1}_{\{  c^{(1)} \neq y_{i+1}\}}     \nonumber  \\
&+ \mathds{1}_{\{c^{(2)} = y_{i+1}\}} - \mathds{1}_{\{c^{(2)} \neq y_{i+1} \}} \big{)} \Big\} +  \beta n^2  \sum_{(x, y), \gamma_{xy} \ni e}   \pi(x, c^{(2)}, y)   \nonumber \\
&\times \exp \Big\{  \tfrac{1}{T} \big{(} -\mathds{1}_{\{x_{i-1} = c^{(2)}\}} + \mathds{1}_{\{x_{i-1} \neq c^{(2)}\}}  + \mathds{1}_{\{x_{i-1} = c^{(1)}\}} - \mathds{1}_{\{x_{i-1} \neq c^{(1)}\}}    \big{)} \Big\} \nonumber\\ 
&=  \alpha n^2 \sum_{w \in \chi:\,w_i=c^{(1)}}  \pi(w) \exp \Big\{ \tfrac{1}{T} \big{(} -\mathds{1}_{\{c^{(1)}=w_{i+1}\}} + \mathds{1}_{\{  c^{(1)} \neq w_{i+1}\}}  + \mathds{1}_{\{c^{(2)} = w_{i+1}\}}      \nonumber   \\
& - \mathds{1}_{\{c^{(2)} \neq w_{i+1} \}}  \big{)}  \Big\}          +  \beta n^2 \sum_{w \in \chi:\, w_i=c^{(2)} }   \pi(w) \exp  \Big\{\tfrac{1}{T} \big{(} -\mathds{1}_{\{w_{i-1} = c^{(2)}\}}  \nonumber  \\
&+ \mathds{1}_{\{w_{i-1} \neq c^{(2)}\}}+ \mathds{1}_{\{w_{i-1} = c^{(1)}\}} - \mathds{1}_{\{w_{i-1} \neq c^{(1)}\}}   \big{)} \Big\} = \alpha n^2 A + \beta n^2 B. \nonumber 
\end{align}  
We now turn  to the computation of the two terms $A$ and $B$ on the right side of the previous equation separately: \\
In order to compute $A$ we generalize some symmetry argument from Shiu and Chen (see \cite{Shiu chen})  to the three state case. In this situation we define three spaces $W^{(k)}=\{ w \in \chi,\,\, w_i=c^{(1)},\,\, w_{i+1}=c^{(k)}\}$ for $k \in \{1, 2, 3 \}$. In order to compute their $\pi$ measure  we will establish some equations between those numbers $\pi(W^{(1)}), \pi(W^{(2)})$ and $\pi(W^{(3)})$.\\
We now establish some identification between the elements from $W^{(1)}$ and the elements of $W^{(2)}$ respective $W^{(3)}$.\\
For any vertex $\xi^{1}\in W^{(1)}$, there exist  a unique vertex $\xi^{2} \in W^{(2)}$ such that:
\begin{itemize}
\item  If $k <i$ , $\xi^{1}_k = \xi^{2}_k$.
\item  If $k>i+1$, then:
\begin{itemize}
\item  If  $\xi^{1}_k = c^{(1)}$ then $\xi^{2}_k = c^{(2)}$.
\item  If $\xi^{1}_k=c^{(2)}$ then  $\xi^{2}_k = c^{(1)}$.
\item If $\xi^{1}_k = c^{(3)}$  then $ \xi^{2}_k = c^{(3)}$.
\end{itemize}
\end{itemize}
Similarly, for any vertex $\xi^1\in W^{(1)}$, there exist a unique vertex $\xi^3 \in W^{(3)}$ such that :
\begin{itemize}
\item  If $k <i$ , $\xi^1_k = \xi^3_k$.
\item  If $k>i+1$, then:
\begin{itemize}
\item  If  $\xi^1_k = c^{(1)}$ then $\xi^3_k = c^{(3)}$.
\item  If $\xi^1_k=c^{(3)}$ then  $\xi^3_k = c^{(1)}$.
\item If $\xi^1_k = c^{(2)}$ then $ \xi^2_k = c^{(2)}$.
\end{itemize}
\end{itemize}
Those relations yield that  \,$\pi(\xi^1)= e^{\frac{2}{T}} \pi(\xi^2)= e^{\frac{2}{T}} \pi(\xi^3).$
Therefore, we obtain: 
\begin{align}
 \displaystyle \sum_{w \in W^{(1)}} \pi(w) = e^{\frac{2}{T}} \displaystyle \sum_{w \in W^{(2)}} \pi(w) = e^{\frac{2}{T}} \displaystyle \sum_{w \in W^{(3)}} \pi(w). \label{form1}
\end{align}
On the other hand  we have also: 
\begin{align}
\sum_{w \in W^{(1)}} \pi(w)+ \sum_{w \in W^{(2)} } \pi(w) + \sum_{w \in W^{(3)}} \pi(w) = \frac{1}{3}. \label{form2}
\end{align}  
From equations (\ref{form1}) and (\ref{form2}) we deduce:$$\sum_{w \in W^{(1)}} \pi(w)= \frac{1}{3(1 +2e^{\frac{-2}{T}})}\,\, \mbox{,}\,\,\, \sum_{w \in W^{(2)} } \pi(w) =  \frac{e^{\frac{-2}{T}}}{3(1 + 2e^{\frac{-2}{T}} ) }\,\,\,\,\,\mbox{and}\,\,\,\,\, \sum_{w \in W^{(3)}} \pi(w) =  \frac{e^{\frac{-2}{T}}}{3(1 + 2e^{\frac{-2}{T}} ) }.$$
The subdivision of the sum in $A$ to three sums over the sets $W^{(1)}$, $W^{(2)}$ and $W^{(3)}$ gives
\begin{align}
 A=\displaystyle\sum_{w \in \chi:\, w_i=c^{(1)}}\pi(w) \exp \Big\{ \tfrac{1}{T} \big{(}  \mathds{1}_{\{c^{(2)} = w_{i+1}\}} - \mathds{1}_{\{c^{(2)} \neq w_{i+1} \}} -\mathds{1}_{\{c^{(1)}=w_{i+1}\}} + \mathds{1}_{\{  c^{(1)} \neq w_{i+1}\}}   \big{)} \Big\} &= \frac{1}{3}. \label{equation-p1}
\end{align}
With the same tricks, we obtain the same result for $B$
\begin{align}
B=\displaystyle \sum_{w \in \chi:\, w_i=c^{(2)}}\pi(w) \exp \Big\{ \tfrac{1}{T} \big{(}\mathds{1}_{\{w_{i-1} = c^{(1)}\}} - \mathds{1}_{\{w_{i-1} \neq c^{(1)}\}} -\mathds{1}_{\{w_{i-1} = c^{(2)}\}} + \mathds{1}_{\{w_{i-1} \neq c^{(2)}\}}   \big{)} \Big\}&= \frac{1 }{3}. \label{equation-p2}
\end{align}
Using equations  $(\ref{equation-p1})$ and $(\ref{equation-p2})$ we obtain:
\begin{align}
Q(e)^{-1} \sum_{(x, y):\, \gamma_{xy} \ni e} |\gamma_{xy}| \pi(x) \pi(y) &\leq \alpha n^2 A + \beta n^2 B  = \frac{n^2}{3}(\alpha+\beta). \nonumber 
\end {align}
The sites $z_{i-1}$ and $z_{i+1}$ take the values $c^{(1)}, c^{(2)}$ or $c^{(3)}$.  The worst value of $\alpha+\beta$ is obtained when $z_{i-1}=c^{(3)}=z_{i+1}$ . In this case we have: 
\begin{align}
Q(e)^{-1} \sum_{(x, y):\, \gamma_{xy} \ni e} | \gamma_{xy}| \pi(x) \pi(y) &\leq \frac{n^2}{3} (2+e^{\frac{4}{T} } ) . \label{eqqq1}
\end {align}
b) If $i=1$ then the configurations $x$ and $e^-$ coincide, from equations (\ref{kappa})  and  (\ref{transition1}) we obtain:
\begin{align}
\kappa &= \max_{e\in E}\,\, Q(e)^{-1} \sum_{\gamma_{xy} \ni e} | \gamma_{xy} |  \pi(x) \pi(y) \leq  n \sum_{\gamma_{xy} \ni e} \frac{\pi(x) \pi(y)}{\frac{1}{n} \pi(e^-) P(e^-, e^+)}   \nonumber \\
&\leq n^2(1+ 2e^{\frac{2}{T}}) \sum_{y \in \chi:\, y_1=c^{(2)}} \pi(y) =\frac{n^2}{3}(2+ e^{\frac{2}{T}}). \label{eqqq2}
\end{align}
For $i=n$, then the configurations  $y$ and $e^+$ coincide and some computation gives a similar result as in equation (\ref{eqqq2}) .\\
By regrouping the results in equations (\ref{eqqq1}) and (\ref{eqqq2}) we obtain an upper bound for the constant $\kappa$ defined in equation (\ref{kappa}) as follows:
\begin{align}
\kappa &= \max_{e\in E} Q(e)^{-1} \sum_{\gamma_{xy} \ni e} | \gamma_{xy} |  \pi(x) \pi(y) \leq  \frac{n^2}{3} (2+e^{\frac{4}{T} } ). \label{pre-res}
\end{align} 
\begin{rem}
The above computation was done  for the situation where $x_i=c^{(1)}$ and $y_i=c^{(2)}$. Obviously we obtain the same result in the other  cases, ie.: $x_i=c^{(1)}$ and $y_i=c^{(3)}$, etc $\cdots$.
\end{rem}
Finally, from  inequalities (\ref{Geometric bound}) and (\ref{pre-res}) we obtain an upper bound for $\beta_1$ which finishes the proof.
\subsection{Proof of theorem $3$}
We follow the same approach as in the case where $\chi=\{c^{(1)}, c^{(2)}, c^{(3)}\}^n$. We define an edge $e=(e^-, e^+)$ as in section $1.2$  and then distinguish two cases:\\
a) For $i \neq \{1, n\}$, we pass in equations  (\ref{transition1}), (\ref{notation}) and (\ref{eq11}) from the case of three colors where $\chi=\{c^{(1)}, c^{(2)}, c^{(3)}\}^n$ to the case of multiple colors where $\chi=\{c^{(1)}, \cdots, c^{(N)}\}^n$. This yields:
\begin{align}
Q(e)^{-1} \sum_{(x, y):\, \gamma_{xy} \ni e} |  \gamma_{xy} | \pi(x) \pi(y) &\leq  \alpha n^2 \sum_{w \in \chi:\, w_i=c^{(1)}}  \pi(w) \exp \Big\{ \tfrac{1}{T} \big{(} -\mathds{1}_{\{c^{(1)}=w_{i+1}\}} + \mathds{1}_{\{  c^{(1)} \neq w_{i+1}\}} +  \mathds{1}_{\{c^{(2)} = w_{i+1}\}   } \nonumber   \\ 
& - \mathds{1}_{\{c^{(2)} \neq w_{i+1} \}}  \big{)}  \Big\}          +  \beta n^2 \sum_{w \in \chi:\, w_i=c^{(2)}}   \pi(w) \exp  \Big\{\tfrac{1}{T} \big{(} -\mathds{1}_{\{w_{i-1} = c^{(2)}\}}  \nonumber   \\
&+ \mathds{1}_{\{w_{i-1} \neq c^{(2)}\}}+ \mathds{1}_{\{w_{i-1} = c^{(1)}\}} - \mathds{1}_{\{w_{i-1} \neq c^{(1)}\}}   \big{)} \Big\}= \alpha n^2 A' + \beta n^2 B'. \nonumber
\end{align}
where
\begin{enumerate}[label= \roman*)]
\item  $\alpha = \exp{ \Big\{  \tfrac{1}{T} \big{(} \mathds{1}_{\{z_{i-1} = c^{(2)}\}}  - \mathds{1}_{\{z_{i-1} \neq c^{(2)}\}}  } - \mathds{1}_{\{z_{i-1} = c^{(1)}\}} + \mathds{1}_{\{z_{i-1} \neq c^{(1)}\}} \}  \big{)} \Big\},$ \\
\item  $\beta=\exp{ \Big\{  \tfrac{1}{T} \big{(} -\mathds{1}_{\{z_{i-1} = c^{(1)}\}} + \mathds{1}_{\{z_{i-1} \neq c^{(1)}\}  } - \mathds{1}_{\{c^{(2)} = z_{i+1}\}} + \mathds{1}_{\{c^{(2)} \neq z_{i+1}\}} \big{)} \Big\} }\\  \times \displaystyle \sum \limits_{\underset{j \neq 2}{j=1}}^N \exp \Big\{ \tfrac{1}{T} \big{(} \mathds{1}_{\{z_{i-1} = c^{(j)}\}} - \mathds{1}_{\{z_{i-1} \neq c^{(j)}\}}+ \mathds{1}_{\{c^{(j)} = z_{i+1}\}} - \mathds{1}_{\{c^{(j)} \neq z_{i+1}\}} \big{)}\Big\}.$
\end{enumerate}
To compute the term $A'$ we define for $k=\{1, \cdots, N\}$ the spaces $W^{(k)}=\{w\in\chi: w_i=c^{(1)},\,w_{i+1}=c^{(k)}\}$  and we consider some symmetry arguments as above:\\
For any vertex $\xi^1\in W^{(1)}$ there exist a unique $\xi^l \in W^{(l)}$ where $l \in \{2, \cdots, N\}$ such that:\begin{itemize}
\item  If $k <i$ , $\xi^1_k = \xi^l_k$.
\item  If $k>i+1$, then:
\begin{itemize}
\item  If  $\xi^1_k = c^{(1)}$ then $\xi^l_k = c^{(l)}$.
\item If $\xi^1_k=c^{(l)}$ then  $\xi^l_k = c^{(1)}$.
\item If $\xi^1_k = c^{(\tilde{l})}$ where $\tilde{l}\neq \{1, l\}$ then $ \xi^l_k =c^{(\tilde{l})}$.
\end{itemize}
\end{itemize}
Then equations (\ref{form1}) and (\ref{form2})  becomes in the case of $N$ colors:
\begin{align}
 \displaystyle \sum_{w \in W^{(1)}} \pi(w) = e^{\frac{2}{T}} \displaystyle \sum_{w \in W^{(2)}} \pi(w) = \cdots =e^{\frac{2}{T}} \displaystyle \sum_{w \in W^{(N)}} \pi(w) , \label{form1-1}
\end{align}
\begin{align}
\sum_{w \in W^{(1)} } \pi(w)+ \cdots + \sum_{w \in W^{(N)}} \pi(w) = \frac{1}{N}. \label{form2-1}
\end{align}  
From equations (\ref{form1-1}) and (\ref{form2-1}) and with the same tricks used to obtain equation (\ref{equation-p1}) we get:
\begin{align}
A'= \displaystyle\sum_{w \in \chi:\,w_i=c^{(1)}}  \pi(w) \exp \Big\{ \tfrac{1}{T} \big{(}\mathds{1}_{\{c^{(2)} = w_{i+1}\}} - \mathds{1}_{\{c^{(2)} \neq w_{i+1} \}} -\mathds{1}_{\{c^{(1)}=w_{i+1}\}} + \mathds{1}_{\{  c^{(1)} \neq w_{i+1}\}}  \big{)} \Big\} &= \frac{1}{N}. \label{equation-p11}
\end{align}
With the same tricks, we obtain the same result for $B'$ 
\begin{align}
B'=\displaystyle \sum_{w \in \chi:\, w_i=c^{(2)}}   \pi(w) \exp \Big\{ \tfrac{1}{T} \big{(}\mathds{1}_{\{w_{i-1} = c^{(1)}\}} - \mathds{1}_{\{w_{i-1} \neq c^{(1)}\}} -\mathds{1}_{\{w_{i-1} = c^{(2)}\}} + \mathds{1}_{\{w_{i-1} \neq c^{(2)}\}}   \big{)} \Big\}&= \frac{1 }{N}. \label{equation-p22}
\end{align}
The application of equations (\ref{equation-p11}) and (\ref{equation-p22}) gives:
\begin{align}
Q(e)^{-1} \sum_{(x, y):\, \gamma_{xy} \ni e} |  \gamma_{xy} | \pi(x) \pi(y) &\leq \alpha n^2 A' + \beta n^2 B'  = \frac{n^2}{N}(\alpha + \beta).\nonumber
\end{align}
The worst value of $\alpha+\beta$ is obtained when $z_{i-1}=c^{(l)} =z_{i+1}$  for  $l \in \{3, \cdots, N\}$. In this case we have: 
\begin{align}
Q(e)^{-1} \sum_{(x, y):\, \gamma_{xy} \ni e} | \gamma_{xy}| \pi(x) \pi(y) &\leq \frac{n^2}{N} (N-1+e^{\frac{4}{T} } ) . \label{eqqq111}
\end {align}
b) For the boundary cases, when $i=1$ equation (\ref{eqqq2}) becomes
\begin{align}
 Q(e)^{-1} \sum_{\gamma_{xy} \ni e} | \gamma_{xy} |  \pi(x) \pi(y) \leq \frac{n^2}{N}( N-1 +e^{\frac{2}{T}}). \label{eqqqq2}
\end{align}
Also, we obtain a similar result for the case where $i=n$.\\
Equation (\ref{eqqq111}) and (\ref{eqqqq2}) together  give an upper bound of $\kappa$ defined in (\ref{kappa}) as follow: \begin{align} \kappa = \max_{e \in E} \, Q(e)^{-1}  \sum_{(x, y):\, \gamma_{xy} \ni e} | \gamma_{xy}| \pi(x) \pi(y) \leq \frac{n^2}{N} (N-1 + e^{\frac{4}{T}}).  \label{step} \end{align}
Finally, we apply the upper bound given in  (\ref{step}) in equation (\ref{Geometric bound}) to finish the proof.
\section*{Acknowledgements}
I would like to thank my supervisors  Brice Franke and Mondher Damak  for their help and advice during this work and their availability for answering all my questions.\\
This work was supported by the Tunisian-French cooperation (PHC-UTIQUE) project CMCU2016 Number 16G1505. \\
I also would like to thank the reviewer  for the comments which greatly helped in ameliorating the paper.

\end{document}